\documentclass[12pt]{amsart}
\usepackage[a4paper]{geometry}
 \newtheorem{thm}{Theorem}[section]
 
 \newtheorem{lem}[thm]{Lemma}
 
 \theoremstyle{definition}
 
 \theoremstyle{remark}
 
 \numberwithin{equation}{section}

\begin{document}
\title {\textbf{
Biharmonic hypersurfaces in space forms with three distinct
principal curvatures}}
\author{Ram Shankar Gupta}
\maketitle
\begin{abstract}
In this paper, we have studied biharmonic hypersurfaces in space
form $\overline{M}^{n+1}(c)$ with constant sectional curvature $c$.
We have obtained that biharmonic hypersurfaces $M^{n}$ with at most
three distinct principal curvatures in $\overline{M}^{n+1}(c)$ has
constant mean curvature. We also obtain the full classification of
biharmonic hypersurfaces with at most three distinct principal
curvatures in arbitrary dimension space form
$\overline{M}^{n+1}(c)$.\\
\\
\textbf{AMS2000 MSC Codes:} 53D12, 53C40, 53C42\\
\textbf{Key Words:} Biharmonic submanifolds, Mean curvature vector,
 isoparametric hypersurfaces.
\end{abstract}

\maketitle
\section{\textbf{Introduction}}
 The longstanding well known Chen's conjecture on biharmonic submanifolds states
that a biharmonic submanifold in a Euclidean space is a minimal one
[2]. In particular, Chen proved that there exist no proper
biharmonic surfaces in Euclidean 3-spaces. There are many
non-existence results in Euclidean spaces developed by I. Dimitri´c
in [9, 10]. Later, the Chen's conjecture was verified and found true
for submanifolds of some Euclidean spaces (see [7, 12, 13, 14]).

 In contrast to the submanifolds of Euclidean spaces, Chen's
conjecture is not always true for the submanifolds of the
pseudo-Euclidean spaces (see [3$\sim$6]). However, for hypersurfaces
in pseudo-Euclidean spaces, Chen's conjecture is also right (see [1,
8]).

 For biharmonic hypersurfaces with at most two distinct principal
curvatures the property of having constant mean curvature was proved
in [15] for any space form. This property proved to be the main
ingredient for the following complete classification of proper
biharmonic hypersurfaces with at most two distinct principal
curvatures in the Euclidean sphere.

\textbf{Theorem 1.1 }([15]): Let $M^{m}$ be a proper biharmonic
hypersurface with at most two distinct principal curvatures in
$\mathbb{S}^{m+1}$. Then $M$ is an open part of
$\mathbb{S}^{m}(\frac{1}{\sqrt{2}})$ or of $\mathbb{S}^{m_{1}}
(\frac{1}{\sqrt{2}})\times \mathbb{S}^{m_{2}}
(\frac{1}{\sqrt{2}}), \quad m_{1}+m_{2}=m,\quad m_{1}\neq m_{2}$.\\

\textbf{Proposition 1.2} ([15]): Let $M^{m}$ be a proper biharmonic
hypersurface with constant mean curvature $H$
  in $\mathbb{S}^{m+1}$. Then $M$ has constant scalar curvature, $$s = m^{2}(1 +
k)-2m,$$ where $H^{2}= k$.

For biharmonic hypersurfaces in $4$-dimensional space form the
property of having constant mean curvature was proved in [19] and
the following classification result was obtained
\\ \textbf{Theorem 1.3}([19]): There exist no
compact proper biharmonic hypersurfaces of constant mean curvature
and with three distinct principal curvatures in the unit Euclidean
sphere.\\\textbf{Theorem 1.4}([19]):
 The only
compact proper biharmonic hypersurfaces of $\mathbb{S}^{4}$ are the
hypersphere $\mathbb{S}^{3}(\frac{1}{\sqrt{2}})$ and the torus
$\mathbb{S}^{1} (\frac{1}{\sqrt{2}})\times \mathbb{S}^{2}
(\frac{1}{\sqrt{2}})$. \\

In view of above development, we study the biharmonic hypersurfaces
in $\overline{M}^{n+1}(c)$ with at most three distinct principal
curvatures.

\section{\textbf{Preliminaries}}

   Let ($M^{n}, g$) be a hypersurface isometrically immersed in a $(n+1)$-dimensional space forms
$(\overline{M}^{n+1}(c), \overline g)$ with constant sectional
curvature $c$ and $g = \overline g_{|M}$.

  Let $\overline\nabla $ and $\nabla$ denote linear connections on $\overline{M}^{n+1}(c)$ and $M^{n}$, respectively. Then,
  the Gauss and Weingarten formulae are given by
\begin{equation}
\overline\nabla_{X}Y = \nabla_{X}Y + h(X, Y), \hspace{.2 cm} \forall
\hspace{.2 cm}X, Y \in\Gamma(TM),
\end{equation}
\begin{equation}
\overline\nabla_{X}\xi = -A_{\xi}X,
\end{equation}
where $\xi$ be the unit normal vector to $M$, $h$ is the second
fundamental form and $A$ is the shape operator. It is well known
that the second fundamental form $h$ and shape operator $A$ are
related by
\begin{equation}
\overline{g}(h(X,Y), \xi) = g(A_{\xi}X,Y).
\end{equation}

 The mean curvature vector is given by
\begin{equation}
H = \frac{1}{n} \emph{trace} A.
\end{equation}

 The Gauss and Codazzi equations are given by
\begin{equation}
R(X, Y)Z = c(g(Y, Z) X - g(X, Z) Y)+g(AY, Z) AX - g(AX, Z) AY,
\end{equation}
\begin{equation}
(\nabla_{X}A)Y = (\nabla_{Y}A)X,
\end{equation}
respectively, where $R$ is the curvature tensor and
\begin{equation}
(\nabla_{X}A)Y = \nabla_{X}AY- A(\nabla_{X}Y)
\end{equation}
for all $ X, Y, Z \in \Gamma(TM)$.

 A biharmonic submanifold in a space form $\overline{M}(c)$ is called proper
 biharmonic if it is not minimal. The necessary and sufficient
 conditions for $M$ to be proper biharmonic in $\overline{M}^{n+1}(c)$ [3, 17] is
\begin{equation}
\triangle H -  H (nc-\emph{trace} A^{2}) = 0,
\end{equation}
\begin{equation}
2A \emph{grad} H+ n H \emph{grad} H = 0,
\end{equation}
where $H$ denotes the mean curvature. Also the Laplace operator
$\triangle$ of a scalar valued function $f$ is given by [4]
\begin{equation}
\triangle f = -\sum_{i=1}^{4}(e_{i}e_{i}f - \nabla_{e_{i}}e_{i}f),
\end{equation}
where $\{e_{1}, e_{2}..., e_{n}\}$ is an orthonormal local tangent
frame on $M^{n}$.\\

  We recall that a hypersurface $M^{n}$ in $\mathbb{S}^{n+1}$ is said
to be isoparametric of type $l$ if it has constant principal
curvatures $k_{1}> . . . > k_{l}$ with respective constant
multiplicities $n_{1}, . . . , n_{l},  \quad n = n_{1} + n_{2} + . .
. + n_{l}$. It is known that the number $l$ is either $1, 2, 3, 4$
or $6$. For $l\leq3$, we have the following classification of
compact isoparametric hypersurfaces. If $l = 1$, then $M$ is totally
umbilical. If $l = 2$, then $M = \mathbb{S}^{n_{1}} (r_{1})\times
\mathbb{S}^{n_{2}} (r_{2}), \quad r_{1}^{2} + r_{2}^{2} = 1$ (see
[18]). If $l= 3$, then $n_{1} = n_{2} = n_{3} = 2^{q}, \quad q = 0,
1, 2, 3$ (see
[16]).\\

  Moreover, there exists an angle $\theta$, $0 < \theta < \frac{\pi}{l},$ such that
\begin{equation}k_{\alpha} = \cot( \theta+\frac{(\alpha-1)\pi}{l}, \quad \alpha
= 1,..., l. \end{equation}

 In the next section, we shall need the following result: \\

\textbf{Theorem 2.1} ([11]): A compact hypersurface $M^{m}$ of
constant scalar curvature $s$ and constant mean curvature $H$ in
$\mathbb{S}^{m+1}$ is isoparametric provided it has 3 distinct
principal curvatures everywhere.\\

\section{\textbf{Biharmonic hypersurfaces with three distinct principal curvatures}}

In this section, we study biharmonic hypersurfaces $M$ in space form
$\overline{M}^{n+1}(c)$. We assume that $H$ is not constant. The
hypothesis for $M$ to be proper biharmonic with three distinct
principal curvatures in space form $\overline{M}^{n+1}(c)$ and
non-constant mean curvature, implies the existence of an open
connected subset $U$ of$M$, with grad$_{p}H\neq 0$ for all $p\in U$.
We shall contradict the condition grad$_{p}H\neq 0$, $\forall p\in
U$. From (2.9), it is easy to see that grad$H$ is an eigenvector of
the shape operator $A$ with the corresponding principal curvature $
\frac{-nH}{2}$. We choose $e_{1}$ in the direction of grad$H$ and
therefore shape operator $A$ of hypersurfaces will take the
following form with respect to a suitable frame  $\{e_{1},
e_{2},..., e_{n-1}, e_{n}\}$
\begin{equation} A_{H}= \left(
                            \begin{array}{cccccc}
                              \frac{-nH}{2} & & \\
                              &   \lambda_{2} & \\
                              & &       .. &  \\
                              & & &             .. &\\
                              & & & &           \lambda_{n-1} &\\
                              & & & & &              \lambda_{n} \\
                            \end{array}
                          \right).
\end{equation}
The grad$H$ can be expressed as
\begin{equation}
grad H =\sum_{i=1}^{n} e_{i}(H)e_{i}.
\end{equation}
 As we have taken $e_{1}$ parallel to grad$H$, consequently
\begin{equation}
e_{1}(H)\neq 0, e_{2}(H)= 0, e_{3}(H) = 0,..., e_{n-1}(H)= 0,
e_{n}(H)= 0.
\end{equation}
We express
\begin{equation}
\nabla_{e_{i}}e_{j}=\sum_{k=1}^{n}\omega_{ij}^{k}e_{k}, \hspace{2
cm} i, j = 1, 2, ..., n.
\end{equation}
Using (3.4) and the compatibility conditions
$(\nabla_{e_{k}}g)(e_{i}, e_{i})= 0$ and $(\nabla_{e_{k}}g)(e_{i},
e_{j})= 0$, we obtain

\begin{equation}
\omega_{ki}^{i}=0, \hspace{1 cm} \omega_{ki}^{j}+ \omega_{kj}^{i}
=0,
\end{equation}
for $i \neq j, $ and $i, j, k = 1, 2, ..., n$. \\
\\
Taking $X=e_{i}, Y=e_{j}$ in (2.7) and using (3.1), (3.4), we get
\begin{center}
$(\nabla_{e_{i}}A)e_{j}=e_{i}(\lambda_{j})e_{j}+\sum_{k=1}^{n}\omega_{ij}^{k}e_{k}
(\lambda_{j}-\lambda_{k}).$
\end{center}
Putting the value of $(\nabla_{e_{i}}A)e_{j}$ in (2.6), we find
\begin{center}
$e_{i}(\lambda_{j})e_{j}+\sum_{k=1}^{n}\omega_{ij}^{k}e_{k}
(\lambda_{j}-\lambda_{k})=e_{j}(\lambda_{i})e_{i}+\sum_{k=1}^{n}\omega_{ji}^{k}e_{k}
(\lambda_{i}-\lambda_{k}),$
\end{center}
whereby for $i\neq j=k$ and $i\neq j \neq k$, we obtain
\begin{equation}
e_{i}(\lambda_{j})= (\lambda_{i}-\lambda_{j})\omega_{ji}^{j},
\end{equation}
\begin{equation}
(\lambda_{i}-\lambda_{j})\omega_{ki}^{j}=
(\lambda_{k}-\lambda_{j})\omega_{ik}^{j},
\end{equation}
respectively, for distinct $i, j, k = 1, 2, ..., n.$\\

 Since $\lambda_{1}=\frac{-nH}{2}$, from (3.3), we get

\begin{equation}
e_{1}(\lambda_{1})\neq 0, e_{2}(\lambda_{1})= 0, e_{3}(\lambda_{1})
= 0,..., e_{n-1}(\lambda_{1})= 0, e_{n}(\lambda_{1})= 0.
\end{equation}

Using (3.8), we have
\begin{center}
$[e_{i}, e_{j}](\lambda_{1})= 0, \hspace{2 cm} i, j = 2, ..., n,$
\end{center}
whereby using (3.4), we find
\begin{equation}
\omega_{ij}^{1}= \omega_{ji}^{1},
\end{equation}
for $i \neq j$ and $i, j = 2, ..., n$.

 Now we show that $\lambda_{j}\neq \lambda_{1}, j= 2, 3,..., n.$ In
 fact, if $\lambda_{j}= \lambda_{1}$ for $j\neq 1$, from (3.6), we
 find
\begin{equation}
e_{1}(\lambda_{j})= (\lambda_{1}-\lambda_{j})\omega_{j1}^{j}=0,
\end{equation}
which contradicts the first expression of (3.8).

 Since $M^{n}$ has three distinct principal curvatures, we can
 assume that $
\lambda_{2}= \lambda_{3}=...=\lambda_{n-1}=\lambda \neq
\lambda_{n}.$ From (2.4), we obtain that

 \begin{equation}
\lambda_{n}=\frac {3nH}{2}- (n-2)\lambda, \hspace{1 cm} \lambda \neq
\frac {-nH}{2},\frac {2nH}{n-2},\frac {3nH}{2(n-1)}.
\end{equation}

 Putting $i, j = 2, 3,...,  n-1,$ and $i\neq j$ in (3.6), we get
\begin{equation}
e_{j}(\lambda)= 0, \hspace{1 cm} \quad\mbox{for}\quad j= 2, 3,...,
n-1.
\end{equation}

 Putting $i\neq 1, j = 1$ in (3.6) and using (3.8) and (3.5), we
 find
\begin{equation}
\omega_{1i}^{1}= 0, \hspace{1 cm}  i= 1, 2, 3,..., n.
\end{equation}

 Putting $i = 2, 3,..., n-1, j = n$ in (3.6) and using (3.12), we
 obtain
\begin{equation}
\omega_{ni}^{n}= 0, \hspace{1 cm}  i=  2, 3,..., n-1.
\end{equation}

 Putting $i = 1, j = 2, 3,...,n-1, n,$ in (3.6), we
 have
\begin{equation}
\omega_{n1}^{n}= \frac
{e_{1}(3nH-2(n-2)\lambda)}{-4nH+2(n-2)\lambda},\hspace{.5 cm}
\omega_{j1}^{j}= -\frac {2e_{1}(\lambda)}{nH+2\lambda}, \hspace{1
cm} j= 2, 3,..., n-1.
\end{equation}

 Putting $i = n, j = 2, 3,...,n-1,$ in (3.6), we
 find
\begin{equation}
\omega_{jn}^{j}= \frac {2e_{n}(\lambda)}{3nH-2(n-1)\lambda},
\hspace{1 cm} j= 2, 3,..., n-1.
\end{equation}

 Putting $i = 1, j\neq k, $ and $j, k = 2, 3,...,n-1,$ in (3.7), we
 obtain
\begin{equation}
\omega_{k1}^{j}= 0, \hspace{1 cm} j\neq k, \quad\mbox{and}\quad j, k
= 2, 3,..., n-1.
\end{equation}

 Putting $i = n, j\neq k, $ and $j, k = 2, 3,...,n-1,$ in (3.7), we
 have
\begin{equation}
\omega_{kn}^{j}= 0, \hspace{1 cm} j\neq k, \quad\mbox{and}\quad j, k
= 2, 3,..., n-1.
\end{equation}

 Putting $i = n, j = 1, $ and $ k = 2, 3,...,n-1,$ in (3.7), and using (3.9) we
 get
\begin{equation}
\omega_{kn}^{1}=\omega_{nk}^{1} = 0, \hspace{1 cm} k = 2, 3,...,
n-1.
\end{equation}

 Putting $i = 1, j = n, $ and $ k = 2, 3,...,n-1,$ in (3.7), and
 using (3.9) we  find
\begin{equation}
\omega_{1k}^{n}=\omega_{k1}^{n} = 0, \hspace{1 cm} k = 2, 3,...,
n-1.
\end{equation}

Now, we have the following:
\begin{lem}
 Let $M^{n}$ be an $n$-dimensional biharmonic hypersurface with three distinct principal curvatures and non-constant mean curvature in
 space forms
 $\overline{M}^{n+1}(c)$ , having the shape operator given by (3.1) with respect to suitable orthonormal frame  $\{e_{1},
e_{2},..., e_{n-1}, e_{n}\}$. Then, we obtain

\begin{equation}
 \nabla_{e_{1}}e_{1}= 0, \hspace{.2 cm} \nabla_{e_{i}}e_{1}=
 -\alpha e_{i}, i= 2, 3,..., n-1,\hspace{.2 cm} \nabla_{e_{n}}e_{1}=\beta
 e_{n},
\end{equation}
\begin{equation}
 \nabla_{e_{i}}e_{i}= \alpha e_{1}+ \sum_{i\neq j, j=2}^{n-1}\omega_{ii}^{j}e_{j}-\frac {2e_{n}(\lambda)}{3nH-2(n-1)\lambda}e_{n},  \hspace{.2 cm} i= 2, 3,..., n-1,
\end{equation}
\begin{equation}
 \nabla_{e_{i}}e_{j}= \sum_{i\neq j, k=2}^{n-2}\omega_{ij}^{k}e_{k},\hspace{.2 cm} i, j= 2, 3,..., n-1,
\end{equation}
\begin{equation}
 \nabla_{e_{1}}e_{n}= 0, \hspace{.1 cm} \nabla_{e_{n}}e_{n}= -\beta e_{1}, \hspace{.1 cm} \nabla_{e_{i}}e_{n}= \frac {2e_{n}(\lambda)}{3nH-2(n-1)\lambda}e_{i},
 \hspace{.2 cm} i= 2, 3,..., n-1,
\end{equation}
where $\omega_{ij}^{k}$ satisfies (3.5) for $i, j, k = 1, 2,
3,...,n-1, n,$ and $\alpha=\frac{2e_{1}(\lambda)}{nH+2\lambda},
\beta = \frac{e_{1}(3nH-2(n-2)\lambda)}{-4nH+2(n-2)\lambda}$.
\end{lem}

 Using Lemma 3.1, Gauss equation and comparing the coefficients with
 respect to a orthonormal frame $\{e_{1},
e_{2},..., e_{n-1}, e_{n}\}$, we find the following:

  $\bullet
X=e_{1}, Y=e_{2}, Z=e_{1},$
\begin{equation}
  e_{1}(\alpha)= \alpha^{2}+c -\frac{nH \lambda}{2}.
\end{equation}

  $\bullet
X=e_{1}, Y=e_{2}, Z=e_{n},$
\begin{equation}
 e_{1}(\frac {2e_{n}(\lambda)}{3nH-2(n-1)\lambda})- \alpha \frac {2e_{n}(\lambda)}{3nH-2(n-1)\lambda}=0.
\end{equation}

  $\bullet
X=e_{1}, Y=e_{n}, Z=e_{1},$
\begin{equation}
 e_{1}(\beta) =- \beta^{2}-c +\frac{nH}{2}( \frac {3nH}{2}-(n-2)\lambda).
\end{equation}
  $\bullet X=e_{3}, Y=e_{n}, Z=e_{1},$
\begin{equation}
 e_{n}(\alpha)+ \frac {2e_{n}(\lambda)}{3nH-2(n-1)\lambda}( \alpha + \beta)= 0.
\end{equation}

$\bullet X=e_{n}, Y=e_{2}, Z=e_{n},$
\begin{equation}
 e_{n}(\frac {2e_{n}(\lambda)}{3nH-2(n-1)\lambda})- \alpha \beta -  (\frac {2e_{n}(\lambda)}{3nH-2(n-1)\lambda})^{2}=-c-\lambda (\frac{3nH}{2}-(n-2)\lambda).
\end{equation}
  Using (2.8), (2.10), (3.1) and Lemma 3.1, we find

\begin{equation}
 -e_{1}e_{1}(H)+ \left[(n-2)\alpha -\beta \right]e_{1}(H)+H \left[\frac{n^{2}H^{2}}{4}+ (n-2)\lambda^{2}+(\frac{3nH}{2}-(n-2)\lambda)^{2}\right]-ncH=0.
\end{equation}

 From (3.3) and Lemma 3.1, we obtain
\begin{equation}
 e_{i}e_{1}(H)= 0, \hspace{1 cm} i= 2, 3,..., n-1, n.
\end{equation}
Differentiating $\alpha=\frac{2e_{1}(\lambda)}{nH+2\lambda}, \beta =
\frac{e_{1}(3nH-2(n-2)\lambda)}{-4nH+2(n-2)\lambda}$ along $e_{n}$,
we get equations
\begin{center}
$(nH+2\lambda)e_{n}(\alpha)+2\alpha
e_{n}(\lambda)=2e_{n}e_{1}(\lambda),$
\end{center}
\begin{center}
$(-4nH+2(n-2)\lambda)e_{n}(\beta)=-2(n-2)e_{n}e_{1}(\lambda)-2(n-2)\beta
e_{n}(\lambda)$
\end{center}
respectively and eliminating $e_{n}e_{1}(\lambda)$, we have
\begin{center}
$(-4nH+2(n-2)\lambda)e_{n}(\beta)=-
(n-2)(nH+2\lambda)e_{n}(\alpha)-2(n-2)(\alpha+\beta)e_{n}(\lambda).$
\end{center}

 Putting the value of $e_{n}(\alpha)$ from (3.28) in the
above equation, we find
\begin{center}
 $e_{n}(\beta)= \frac{4e_{n}(\lambda)n(n-2)(\alpha+\beta)(\lambda-H)}{(-4nH+2(n-2)\lambda)(3nH-(2n-2)\lambda)}$.
\end{center}
Differentiating (3.30) along $e_{n}$ and using (3.31), (3.28) and
$e_{n}(\beta)$, we get
\begin{equation}
 e_{n}(\lambda)\left[\frac{4(\alpha+\beta)e_{1}(H)}{-4nH+2(n-2)\lambda}+H((2n-2)\lambda-3nH)\right]= 0.
\end{equation}

We claim that $e_{n}(\lambda)= 0$. Indeed, if $e_{n}(\lambda)\neq
0$, then

\begin{equation}
 \frac{4(\alpha+\beta)e_{1}(H)}{-4nH+2(n-2)\lambda}+H((2n-2)\lambda-3nH)= 0.
\end{equation}
Now, differentiating (3.33) along $e_{n}$, we have
\begin{equation}
 \frac{8(\alpha+\beta)(nH(14-5n)+4(n-2)(n-1)\lambda)e_{1}(H)}{(-4nH+2(n-2)\lambda)^2 (\frac{3nH}{2}-(n-2)\lambda)}+H((2n-2)= 0.
\end{equation}
Eliminating $e_{1}(H)$ from (3.33) and (3.34), we obtain
\begin{center}
$2(n-1)\lambda-3nH=0$
\end{center}
which is not possible since $ \lambda \neq \frac{3nH}{2(n-1)}$,
consequently, $e_{n}(\lambda)=0$. Therefore, (3.29) reduces to
\begin{equation}
 \alpha \beta = c+\lambda (\frac{3nH}{2}-(n-2)\lambda).
\end{equation}
Now, eliminating $e_{1}e_{1}(H)$ and $e_{1}e_{1}(\lambda)$, using
(3.35), (3.30), (3.27) and (3.25), we obtain
\begin{equation}
 [(10n-2n^{2})\alpha - 4n \beta]e_{1}(H)= \frac
{21n^{3}H^{3}}{2}+6(n^{3}-2n^{2})H\lambda^{2}+(-15n^{3}+18n^{2})H^{2}\lambda-6(n^{2}+n)cH.
\end{equation}
Differentiating (3.36) along $e_{1}$ and using (3.35), (3.30),
(3.27), (3.25) and (3.36), we get
\begin{equation}
\begin{array}{rcl}
 [
(13n^{3}+\frac{11n^{2}}{2})H^{3}+(4n^{3}-14n^{2}+2n+20)H\lambda^{2}+(-15n^{3}+18n^{2}+24n)H^{2}\lambda+cH(2n^{3}\\-16n^{2}-6n)]\alpha
+ [-31n^{2}H^{3}+(-16n^{2}+36n-8)H\lambda^{2}+(42n^{2}-
60n)H^{2}\lambda\\+cH(10n^{2}+6n)] \beta  = e_{1}(H)[\frac
{69n^{2}H^{2}}{2}+(24n -30n^{2})H\lambda
+(6n+4n^{2}-28)\lambda^{2}-c(4n^{2}+20n)].
\end{array}
\end{equation}

Also, we have
\begin{equation}
3ne_{1}(H)= \alpha (n-2)(nH+2\lambda)+ \beta (-4nH+2(n-2)\lambda)
\end{equation}
 Combining (3.37) and (3.38), we obtain

\begin{equation}
\begin{array}{rcl}
 [
(9n^{3}+171n^{2})H^{3}+(16n^{3}+40n^{2}-244n-200)H\lambda^{2}+(-30n^{3}-198n^{2}-516n)H^{2}\lambda\\ -(16n^{2}-8n-160+\frac{224}{n})\lambda^{3}
+cH(20n^{3}-72n^{2}-116n)+c\lambda(16n^{2}+48n-160)]\alpha\\
+
[90n^{2}H^{3}+(56n^{2}-72n-80)H\lambda^{2}+(-126n^{2}+108n)H^{2}\lambda
-(16n^{2}-8n-160+\frac{224}{n})\lambda^{3}\\+cH(28n^{2}-124n)+c\lambda(16n^{2}+48n-160)]
\beta = 0.
\end{array}
\end{equation}
For simplicity, we denote by
\begin{center}
$p_{1}=
(9n^{3}+171n^{2})H^{3}+(16n^{3}+40n^{2}-244n-200)H\lambda^{2}+(-30n^{3}-198n^{2}-516n)H^{2}\lambda
-(16n^{2}-8n-160+\frac{224}{n})\lambda^{3}+cH(20n^{3}-72n^{2}-116n)+c\lambda(16n^{2}+48n-160)$
\end{center}
\begin{center}
$p_{2}=
90n^{2}H^{3}+(56n^{2}-72n-80)H\lambda^{2}+(-126n^{2}+108n)H^{2}\lambda
-(16n^{2}-8n-160+\frac{224}{n})\lambda^{3}+cH(28n^{2}-124n)+c\lambda(16n^{2}+48n-160).$
\end{center}
Therefore, (3.39) can be rewritten as
\begin{equation}
\alpha p_{1} + \beta p_{2}=0.
\end{equation}

 On the other hand, combining (3.38) with
(3.36) and using (3.35), we find

\begin{equation}
\alpha^{2} (n-2)(10-2n)(nH+2\lambda)- 4 \beta^{2}
(-4nH+2(n-2)\lambda)= L,
\end{equation}
where $L$ is given by
\begin{center}
$L=
\frac{63n^{3}H^{3}}{2}+(28n^{3}-106n^{2}+100n)H\lambda^{2}+(102n^{2}-51n^{3})H^{2}\lambda
 -(4n^{3}-28n^{2}+64n-48)\lambda^{3}+cH(14n-22n^{2})+c\lambda(4n^{2}-20n+24).$
\end{center}

  Using (3.40) and (3.35), we get
\begin{center}
$\alpha^{2}= -\frac{
p_{2}}{p_{1}}(c+\lambda(\frac{3nH}{2}-(n-2)\lambda)), \hspace{2 cm}
\beta^{2}=
-\frac{p_{1}}{p_{2}}(c+\lambda(\frac{3nH}{2}-(n-2)\lambda))$
\end{center}

  Eliminating $\alpha^{2}$ and $\beta^{2}$ from (3.41), we obtain
\begin{equation}
(c+\frac{3nH\lambda}{2}-(n-2)\lambda^{2})[(14n-2n^{2}-20)(nH+2\lambda)
p_{2}^{2}- 4 p_{1}^{2}(-4nH+2(n-2)\lambda)]= L p_{1}p_{2},
\end{equation}
which is a polynomial equation of degree 9 in terms of $\lambda$ and
$H$.

Now consider an integral curve of $e_{1}$ passing through $p =
\gamma(t_{0})$ as $\gamma(t)$, $t\in I$. Since $e_{i}(H) =
e_{i}(\lambda) = 0$ for $i = 2, . . ., n $ and $e_{1}(H),
e_{1}(\lambda)\neq 0$, we can assume $t = t(\lambda)$ and $H =
H(\lambda)$ in some neighborhood of $\lambda_{0}=\lambda(t_{0})$.
Using (3.38) and (3.40), we have
\begin{equation}\begin{array}{lcl}
\frac{dH}{d\lambda}=\frac{dH}{dt}\frac{dt}{d\lambda}=\frac{e_{1}(H)}{e_{1}(\lambda)}\\
=\frac{2(\alpha (n-2)(nH+2\lambda)+ \beta
(-4nH+2(n-2)\lambda))}{3n\alpha(nH+2\lambda)}\\
=\frac{2(n-2)}{3n}+\frac{p_{1}(4nH-2(n-2)\lambda))}{3n(nH+2\lambda)p_{2}}
\end{array}
\end{equation}
 Differentiating (3.42) with respect to $\lambda$ and substituting $\frac{dH}{d\lambda}$ from (3.43), we get
\begin{equation}\begin{array}{lcl}
f(H, \lambda)=0,
\end{array}
\end{equation}
  another algebraic
equation of degree 12 in terms of $H$ and $\lambda$. We rewrite
(3.42) and (3.44) respectively in the following forms
\begin{equation}\begin{array}{lcl}
\sum_{i=0}^{9}f_{i}(H) \lambda^{i}=0,\quad \sum_{j=0}^{12}g_{j}(H)
\lambda^{j}=0,
\end{array}
\end{equation}
where $f_{i}(H)$ and $g_{j}(H)$ are polynomial functions of $H$. We
eliminate $\lambda^{12}$ between these two polynomials of (3.45) by
multiplying $g_{12}\lambda^{3}$ and $f_{8}$ respectively on the
first and second equations of (3.45), we obtain a new polynomial
equation in $\lambda$ of degree 11. Combining this equation with the
first equation of (3.45), we successively obtain a polynomial
equation in $\lambda$ of degree 10. In a similar way, by using the
first equation of (3.45) and its consequences we are able to
gradually eliminate  $\lambda$. At last, we obtain a non-trivial
algebraic polynomial equation in $H$ with constant coefficients.
Therefore, we conclude that the real function $H$ must be a constant
and we conclude:

 \begin{thm}
 Every  biharmonic  hypersurface $M$ in the space forms
$\overline{M}^{n+1}(c)$ with three distinct principal curvatures
must be of constant mean curvature.
 \end{thm}
 Combining Theorem 3.2 with Theorem 4.1 [15], we obtain that
 \begin{thm}
 Every  biharmonic  hypersurface $M$ in the space forms
$\overline{M}^{n+1}(c)$ with at most three distinct principal
curvatures must be of constant mean curvature.
 \end{thm}
 \begin{thm}
There exist no proper biharmonic hypersurfaces $M$ with at most
three distinct principal curvatures in $\mathbb{H}^{n+1}$ or
$\mathbb{R}^{n+1}$.
 \end{thm}
 $\emph{\textbf{Proof:}}$ Suppose that $M^{n}$ is a
  proper biharmonic hypersurface in $\mathbb{H}^{n+1}$ or $\mathbb{R}^{n+1}$ with at most three distinct principal curvatures. From Theorem
  3.3, we have that mean curvature of $M^{n}$ is constant. From
  (2.8), we get that trace$A^{2}=-n$ or trace$A^{2}=0$, which is not possible and
  proof of the theorem is complete.
 \begin{thm}
 The only compact proper biharmonic hypersurfaces with at most three
 distinct principal curvatures
 of $\mathbb{S}^{n+1}(1)$ are the hypersphere $\mathbb{S}^{n}(\frac{1}{\sqrt{2}})$ and the torus $\mathbb{S}^{n_{1}}(\frac{1}{\sqrt{2}})
 \times \mathbb{S}^{n_{2}}(\frac{1}{\sqrt{2}})$ where $\quad n_{1}+n_{2}=n,\quad
n_{1}\neq n_{2}$.
 \end{thm}
$\emph{\textbf{Proof:}}$ Suppose that $M^{n}$ is a compact proper
biharmonic hypersurface of $\mathbb{S}^{n+1}(1)$ with three distinct
principal curvatures. From Theorem 3.2, we get that $M^{n}$ has
constant mean curvature and, since it satisfies the hypotheses of
Proposition 1.2, we conclude that it also has constant scalar
curvature. We can thus apply Theorem 2.1 and it results that $M^{n}$
is isoparametric in $\mathbb{S}^{n+1}(1)$. From Theorem 1.3, we get
that $M^{n}$ cannot be isoparametric with $l = 3$, and by using
Theorem 1.1 we conclude the proof.

\bibliography{xbib}


Author's address:

\textbf{Ram Shankar Gupta }\\ Associate Professor, Department of
Mathematics, Central University of Jammu, Sainik colony,
Jammu-180011, India.
\\
Assistant Professor, University School of Basic and Applied
Sciences, Guru Gobind Singh Indraprastha University, Sector-16C,
Dwarka, New Delhi-110078,
India.\\
\textbf{Email:} ramshankar.gupta@gmail.com

\end{document}